\documentclass{amsart}

\usepackage{amsmath}
\usepackage{amssymb}

\usepackage{fancyvrb}
\newtheorem{defn}{Definition}
\newtheorem{rem}{Remark}

\usepackage[pdfauthor={Richard J. Mathar},colorlinks=true,citecolor=blue]{hyperref}

\begin{document}
\title{Integrals Associated with the Digamma Integral Representation}

\author{Richard J. Mathar}
\urladdr{http://www.mpia.de/~mathar}
\address{Hoeschstr.\ 7, 52372 Kreuzau, Germany}


\date{\today}

\begin{abstract}
The definite integral with the kernel $x/(x^2+b^2)/[\exp(2\pi x)-1]$ 
integrated from $x=0$ to infinity is the main term 
of a representation of the
Digamma-Function $\psi(b)$, 
the derivative of the logarithm of the Gamma-Function. 
We present relations within the set of integrals 
over $x^n/(x^2+b^2)^j/[\exp(\mu x)-1]^s$ for small 
integer exponents $n$, $j$ and $s$ with the aim
to reduce them all to Polygamma-Functions.
\end{abstract}

\maketitle

\section{Notation}
The integral representation of the digamma function \cite[1.7.2 (27)]{ErdelyiI}
$$
\psi(z) = \log z -\frac{1}{2z}-2\int_0^\infty dt \frac{t}{(t^2+z^2)(e^{2\pi t} -1)}
$$
is one method of computing the integrals \cite[3.415.1]{GR}
\begin{equation}
\int_0^\infty dx \frac{x}{(x^2+b^2)(e^{\mu x} -1)}
=\frac12\left[\log(\frac{b\mu }{2\pi })-\frac{\pi}{b\mu}-\psi(\frac{b\mu}{2\pi})\right]
.
\label{eq.fund}
\end{equation}
Throughout the manuscript $b,\mu >0$.
The theme of this work is the reduction of the following familiy
of integrals
for small positive integer exponents $n$, $j$ and $s$:
\begin{defn} (Associates)
\begin{equation}
I^\pm _{n,j,s}(b,\mu) \equiv \int_0^\infty \frac{x^n}{(x^2+b^2)^j(e^{\mu x}\pm 1)^s}dx
.
\end{equation}
\end{defn}
To some degree one of the two parameters $\mu$ and $b$ is redundant, 
because a trivial substitution 
demonstrates that only their product $b\mu$ is relevant:
\begin{equation}
I^\pm _{n,j,s}(b,\mu) 
=
\mu ^{2j-n-1} I^\pm _{n,j,s}(b\mu,1)
=
b^{n+1-2j} I^\pm _{n,j,s}(1,b\mu)
.
\end{equation}
\begin{defn} (Normal Form)
\begin{equation}
\hat I^\pm _{n,j,s}(b) \equiv \int_0^\infty \frac{x^n}{(x^2+b^2)^j(e^x \pm 1)^s}dx
= I^\pm_{n,j,s}(b,1)
.
\label{def.norm}
\end{equation}
\end{defn}
However, maintaining both parameters $b$ and $\mu$ furnishes
different ways to differentiation ``under the integral'', exploited
later in this manuscript.

\section{Basic Integrals for $I^-$}\label{sec.diffb}

At $j=0$
\cite[3.411.1]{GR}
\begin{equation}
I^-_{n,0,1} = \int_0^\infty \frac{x^n}{e^{\mu x}-1}dx
=\frac{n!}{\mu^{n+1}}\zeta(n+1)
\label{eq.noden}
\end{equation}
in terms of the Riemann $\zeta$-function.

By iterated differentiation of \eqref{eq.fund} with respect to (w.r.t.) $b$ one 
can increase $j$ \cite{GeddesCM}:
\begin{equation}
I^-_{1,2,1}=\int_0^\infty \frac{xdx}{(x^2+b^2)^2 (e^{\mu x}-1)}
= 
-\frac{1}{4b^2}
-\frac{\pi}{4b^3\mu}
+\frac{\mu}{8b\pi}\psi'\left(\frac{b\mu}{2\pi}\right)
. 
\end{equation}
\begin{equation}
I^-_{1,3,1}=\int_0^\infty \frac{xdx}{(x^2+b^2)^3 (e^{\mu x}-1)}
= 
-\frac{1}{8b^4}
-\frac{3\pi}{16b^5\mu}
+\frac{\mu}{32b^3\pi}\psi'\left(\frac{b\mu}{2\pi}\right)
-\frac{\mu^2}{64b^2\pi^2}\psi''\left(\frac{b\mu}{2\pi}\right)
. 
\end{equation}
\begin{multline}
I^-_{1,4,1}=\int_0^\infty \frac{xdx}{(x^2+b^2)^4 (e^{\mu x}-1)}
=
-\frac{1}{12b^6}
-\frac{5\pi}{32b^7\mu}
+\frac{\mu}{64b^5\pi}\psi'\left(\frac{b\mu}{2\pi}\right)
\\
-\frac{\mu^2}{128b^4\pi^2}\psi''\left(\frac{b\mu}{2\pi}\right)
+\frac{\mu^3}{768b^3\pi^3}\psi'''\left(\frac{b\mu}{2\pi}\right)
.
\end{multline}
\begin{multline}
I^-_{1,5,1}=\int_0^\infty \frac{xdx}{(x^2+b^2)^5 (e^{\mu x}-1)}
=
-\frac{1}{16b^8}
-\frac{35\pi}{256b^9\mu}
+\frac{5\mu}{512b^7\pi}\psi'\left(\frac{b\mu}{2\pi}\right)
\\
-\frac{5\mu^2}{1024b^6\pi^2}\psi''\left(\frac{b\mu}{2\pi}\right)
+\frac{\mu^3}{1024b^5\pi^3}\psi'''\left(\frac{b\mu}{2\pi}\right)
-\frac{\mu^4}{12288b^4\pi^4}\psi^{(4)}\left(\frac{b\mu}{2\pi}\right)
.
\end{multline}
\begin{multline}
I^-_{1,6,1}=\int_0^\infty \frac{xdx}{(x^2+b^2)^6 (e^{\mu x}-1)}
=
-\frac{1}{20b^{10}}
-\frac{63\pi}{512b^{11}\mu}
+\frac{7\mu}{1024b^9\pi}\psi'\left(\frac{b\mu}{2\pi}\right)
\\
-\frac{7\mu^2}{2048b^8\pi^2}\psi''\left(\frac{b\mu}{2\pi}\right)
+\frac{3\mu^3}{4096b^7\pi^3}\psi'''\left(\frac{b\mu}{2\pi}\right)
\\
-\frac{\mu^4}{12288b^6\pi^4}\psi^{(4)}\left(\frac{b\mu}{2\pi}\right)
+\frac{\mu^5}{245760b^5\pi^5}\psi^{(5)}\left(\frac{b\mu}{2\pi}\right)
.
\end{multline}
In summary: If expansion coefficients $\gamma^-_{j,i}$ for the normalized integrals are defined as
\begin{equation}
\hat I^-_{1,j,1}
\equiv
-\gamma^-_{j,-1}\frac{\pi}{(b\mu)^{2j-1}}
-\gamma^-_{j,0}\frac{1}{(b\mu)^{2j-2}}
-\sum_{i=1}^{j-1} (-)^i\gamma^-_{j,i}\frac{1}{(b\mu)^{2j-i-2}\pi^i}\psi^{(i)}(\frac{b\mu}{2\pi})
,
\end{equation}
they are recursively constructable
for $j>2$ by
\begin{equation}
\gamma^-_{j,i} = \left\{
\begin{array}{rr}
\frac{2j-3}{2j-2}\gamma^-_{j-1,i}, & i=-1;\\
\frac{j-2}{j-1}\gamma^-_{j-1,i}, & i=0;\\
\frac{2j-i-4}{2j-2}\gamma^-_{j-1,i},& i=1;\\
\frac{2j-i-4}{2j-2}\gamma^-_{j-1,i} + \frac{1}{4(j-1)}\gamma^-_{j-1,i-1},& 2\le i < j-1;\\
\frac{1}{4(j-1)}\gamma^-_{j-1,i-1},& i = j-1.
\end{array}
\right.
\end{equation}
starting
at $\gamma^-_{2,-1}=\gamma^-_{2,0}=1/4$, $\gamma^-_{2,1}=1/8$.

\begin{rem}
Numerical aspects of evaluating the polygamma functions
are not discussed in this paper. Most algorithms use the equations
derived from the functional equation or the duplication formula of the $\Gamma$-function
to reduce the arguments
to the unit interval
\cite{WimpMathComp15, Cody1973, ReinartzArxiv1605,CausleyNA90,ChoiJPA40}.
\end{rem}
\begin{rem}
The integral kernels are functions of $b^2$, even functions of $b$, whereas the right hand sides
reveal no such symmetry. These are hidden by the analogues of the reflection formula of the $\Gamma$-function,
\begin{equation}
\psi'(-x) +\psi'(x)
= \frac{\pi^2}{\sin^2(\pi x)}+\frac{1}{x^2}
\end{equation}
and its higher derivatives.
This does not help to simplify the notation here.
\end{rem}

\section{Basic Integrals for $I^+$} \label{sec.1j1}
The variant with $+1$ in the denominator is reduced to the variant with $-1$ for arguments $\mu$ and $2\mu$
by factorizing $e^{\mu x}-1=(e^{\mu x/2}+1)(e^{\mu x/2}-1)$ followed by partial fraction decomposition:
\begin{equation}
\int_0^\infty \frac{x^n dx}{(x^2+b^2)^j (e^{\mu x}+1)}
=
\int_0^\infty \frac{x^n dx}{(x^2+b^2)^j (e^{\mu x}-1)}
-2\int_0^\infty \frac{x^n dx}{(x^2+b^2)^j (e^{2\mu x}-1)}
;
\end{equation}
\begin{equation}
I^+_{n,j,1}(b,\mu)
=
I^-_{n,j,1}(b,\mu)
-2I^-_{n,j,1}(b,2\mu).
\end{equation}
More compact polygamma expansions are obtained starting from \cite[3.415.3]{GR}
\begin{equation}
\int_0^\infty \frac{x}{(x^2+b^2)(e^{\mu x}+1)}dx
=\frac12\left[\psi(\frac{b\mu}{2\pi}+\frac12)-\ln(\frac{b\mu}{2\pi})\right] 
\end{equation}
and iterated differentiation w.r.t.\ b which raises the parameter $j$:
\begin{equation}
I^+_{1,2,1}=\int_0^\infty \frac{xdx}{(x^2+b^2)^2 (e^{\mu x}+1)}
= 
\frac{1}{4b^2}
-\frac{\mu}{8b\pi}\psi'\left(\frac{b\mu}{2\pi}+\frac12\right)
. 
\end{equation}
\begin{equation}
I^+_{1,3,1}=\int_0^\infty \frac{xdx}{(x^2+b^2)^3 (e^{\mu x}+1)}
= 
\frac{1}{8b^4}
-\frac{\mu}{32b^3\pi}\psi'\left(\frac{b\mu}{2\pi}+\frac12\right)
+\frac{\mu^2}{64b^2\pi^2}\psi''\left(\frac{b\mu}{2\pi}+\frac12\right)
. 
\end{equation}
\begin{multline}
I^+_{1,4,1}=\int_0^\infty \frac{xdx}{(x^2+b^2)^4 (e^{\mu x}+1)}
= 
\frac{1}{12b^6}
-\frac{\mu}{64b^5\pi}\psi'\left(\frac{b\mu}{2\pi}+\frac12\right)
\\
+\frac{\mu^2}{128b^4\pi^2}\psi''\left(\frac{b\mu}{2\pi}+\frac12\right)
-\frac{\mu^3}{768b^3\pi^3}\psi'''\left(\frac{b\mu}{2\pi}+\frac12\right)
. 
\end{multline}
\begin{multline}
I^+_{1,5,1}=\int_0^\infty \frac{xdx}{(x^2+b^2)^5 (e^{\mu x}+1)}
= 
\frac{1}{16b^8}
-\frac{5\mu}{512b^7\pi}\psi'\left(\frac{b\mu}{2\pi}+\frac12\right)
\\
+\frac{5\mu^2}{1024b^6\pi^2}\psi''\left(\frac{b\mu}{2\pi}+\frac12\right)
-\frac{\mu^3}{1024b^5\pi^3}\psi'''\left(\frac{b\mu}{2\pi}+\frac12\right)
+\frac{\mu^4}{12288b^4\pi^4}\psi^{(4)}\left(\frac{b\mu}{2\pi}+\frac12\right)
. 
\end{multline}
\begin{multline}
I^+_{1,6,1}=\int_0^\infty \frac{xdx}{(x^2+b^2)^6 (e^{\mu x}+1)}
= 
\frac{1}{20b^{10}}
-\frac{7\mu}{1024 b^9\pi}\psi'\left(\frac{b\mu}{2\pi}+\frac12\right)
\\
+\frac{7\mu^2}{2048b^8\pi^2}\psi''\left(\frac{b\mu}{2\pi}+\frac12\right)
-\frac{3\mu^3}{4096b^7\pi^3}\psi'''\left(\frac{b\mu}{2\pi}+\frac12\right)
\\
+\frac{\mu^4}{12288b^6\pi^4}\psi^{(4)}\left(\frac{b\mu}{2\pi}+\frac12\right)
-\frac{\mu^5}{245760 b^5\pi^4}\psi^{(5)}\left(\frac{b\mu}{2\pi}+\frac12\right)
. 
\end{multline}

In overview: If expansion coefficients $\gamma^+_{j,i}$ for the normalized integrals are defined as
\begin{equation}
\hat I^+_{1,j,1}
\equiv
\gamma^+_{j,0}\frac{1}{(b\mu)^{2j-2}}
+\sum_{i=1}^{j-1} (-)^i\gamma^+_{j,i}\frac{1}{(b\mu)^{2j-i-2}\pi^i}\psi^{(i)}(\frac{b\mu}{2\pi}+\frac12)
,
\end{equation}
they are recursively constructable
for $j>2$ by
\begin{equation}
\gamma^+_{j,i} = \left\{
\begin{array}{rr}
\frac{j-2}{j-1}\gamma^+_{j-1,i}, & i=0;\\
\frac{2j-i-4}{2j-2}\gamma^+_{j-1,i},& i=1;\\
\frac{2j-i-4}{2j-2}\gamma^+_{j-1,i} + \frac{1}{4(j-1)}\gamma^+_{j-1,i-1},& 2\le i < j-1;\\
\frac{1}{4(j-1)}\gamma^+_{j-1,i-1},& i = j-1.
\end{array}
\right.
\end{equation}
starting
at $\gamma^+_{2,0}=1/4$, $\gamma^+_{2,1}=1/8$.
Obviously $\gamma^+_{j,i}=\gamma^-_{j,i}$ for $j\ge 2$, $i\ge 0$.

\section{Exponent Shifts with Partial Fractions}\label{sec.binom}
High integer exponents $n$ in the numerator are lowered by partial fractions
in steps of two until they reach $s$ or $s+1$:
\begin{equation}
I^\pm _{n,j,s}
=
I^\pm _{n-2,j-1,s}
-b^2 I^\pm _{n-2,j,s},\quad n\ge 2+s
.
\label{eq.lown}
\end{equation}
Complementary, for negative $j$ the exponent is increased via
\begin{equation}
I^\pm _{n,-j,s}=\int_0^\infty \frac{x^n(x^2+b^2)^j}{(e^{\mu x}\pm 1)^s} dx
=
\sum_{l=0}^j \binom{j}{l}
b^{2(n-l)} I^\pm _{n+2l,0,s}.
\end{equation}
(The $I^+_{n,-j,s}$ needs $n+2j\ge s$ for $b=0$ and
$n\ge s$ for $b\neq 0$ to converge.)
For $s=1$ these are evaluated via \eqref{eq.noden}.

\section{Derivatives w.r.t.\ $\mu$}\label{sec.dmu}
The derivative of \eqref{eq.noden} w.r.t.\ $\mu$ is
\begin{equation}
I_{n+1,0,1}+I_{n+1,0,2}
=
\frac{(n+1)!}{\mu ^{n+2}}\zeta(n+1)
.
\end{equation}
The derivative of \eqref{eq.fund} w.r.t.\ $\mu$ is
\begin{equation}
-\int_0^\infty dx \frac{x^2 e^{\mu x}}{(x^2+b^2)(e^{\mu x} -1)^2}
=\frac{1}{2\mu}+\frac{\pi}{2b\mu^2}-\frac{b}{4\pi}\psi'(\frac{b\mu}{2\pi})
.
\label{eq.mu1}
\end{equation}
The powers of $e^{\mu x}$ that appear in the numerators of these derivatives
are eliminated by partial fraction decompositions along the scheme
\begin{equation}
\frac{d}{d\mu} \frac{1}{e^{\mu x}-1}
= -x[\frac{1}{e^{\mu x}-1}+\frac{1}{(e^{\mu x}-1)^2}]
;
\end{equation}
\begin{equation}
\frac{d}{d\mu^2} \frac{1}{e^{\mu x}-1}
= x^2[\frac{1}{e^{\mu x}-1}
+\frac{3}{(e^{\mu x}-1)^2}
+\frac{2}{(e^{\mu x}-1)^3}
]
;
\end{equation}
\begin{equation}
\frac{d}{d\mu^3} \frac{1}{e^{\mu x}-1}
= -x^3[\frac{1}{e^{\mu x}-1}
+\frac{7}{(e^{\mu x}-1)^2}
+\frac{12}{(e^{\mu x}-1)^3}
+\frac{6}{(e^{\mu x}-1)^4}
]
;
\end{equation}
\begin{table}
\begin{tabular}{r|rrrrrrrr}
$l\backslash k$ & 1 & 2 &3 & 4 & 5 & 6 & 7\\
\hline
1 & 1 \\
2 & 1 & 1\\
3 & 1 & 3 & 2\\
4 & 1& 7 & 12 & 6\\
5 & 1 & 15 & 50 & 60 & 24 \\
6 & 1 & 31 & 180 & 390 & 360 & 120 \\
7 & 1 & 63 & 602 & 2100 & 3360 & 2520 & 720 \\
\end{tabular}
\caption{Expansion coefficients $\alpha_{l,k}$ \cite[A028246]{sloane}.
Stirling Numbers of the Second Kind multiplied by $(k-1)!$.}
\label{tab.alpha}
\end{table}
\begin{eqnarray}
\frac{d^l}{d\mu^l} \frac{1}{e^{\mu x}-1}
&=& (-x)^l\sum_{k=1}^{l+1} \frac{\alpha_{l+1,k}}{(e^{\mu x}-1)^k};
\\
\frac{d^l}{d\mu^l} \frac{1}{e^{\mu x}+1}
&=& x^l\sum_{k=1}^{l+1} (-)^{l+k+1}\frac{\alpha_{l+1,k}}{(e^{\mu x}+1)^k},
\label{eq.powexpmu}
\end{eqnarray}
with positive integer coefficients (Table \ref{tab.alpha}) \cite{SofoSci17}
\begin{equation}
\alpha_{l,k} \equiv \frac{1}{k} \sum_{i=0}^k (-)^{k-i}\binom{k}{i}i^l.
\end{equation}
This reduces the $\mu$-derivatives into the name space of the $I^\pm$.
\eqref{eq.mu1} is rewritten as
\begin{equation}
-\alpha_{2,1}I^-_{2,1,1}
-\alpha_{2,2}I^-_{2,1,2}
=\frac{1}{2\mu}+\frac{\pi}{2b\mu^2}-\frac{b}{4\pi}\psi'(\frac{b\mu}{2\pi})
,
\label{eq.dmu2}
\end{equation}
and higher derivatives $\partial/\partial \mu$ are
\begin{equation}
\alpha_{3,1}I^-_{3,1,1}
+\alpha_{3,2}I^-_{3,1,2}
+\alpha_{3,3}I^-_{3,1,3}
=-\frac{1}{2\mu^2}-\frac{\pi}{b\mu^3}-\frac{b^2}{8\pi^2}\psi''(\frac{b\mu}{2\pi})
;
\label{eq.dmu3}
\end{equation}
\begin{equation}
-\alpha_{4,1}I^-_{4,1,1}
-\alpha_{4,2}I^-_{4,1,2}
-\alpha_{4,3}I^-_{4,1,3}
-\alpha_{4,4}I^-_{4,1,4}
=\frac{1}{\mu^3}+\frac{3\pi}{b\mu^4}-\frac{b^3}{16\pi^3}\psi'''(\frac{b\mu}{2\pi})
;
\end{equation}
\begin{equation}
\sum_{s=1}^n \alpha_{n,s}I^-_{n,1,s}
=- \frac{(n-2)!}{2\mu^{n-1}}-\frac{\pi}{2b} \frac{(n-1)!}{\mu^n}+(-)^n\frac12(\frac{b}{2\pi})^{n-1}\psi^{(n-1)}(\frac{b\mu}{2\pi})
.
\label{eq.sumalpha}
\end{equation}
The same in terms of the normalized integrals:
\begin{equation}
\sum_{s=1}^n \alpha_{n,s}\hat I^-_{n,1,s}
=- \frac{(n-2)!}{2}-\frac{\pi(n-1)!}{2b\mu} +(-)^n\frac12(\frac{b\mu}{2\pi})^{n-1}\psi^{(n-1)}(\frac{b\mu}{2\pi})
.
\end{equation}
This sum rule is not sufficient to bootstrap the integrals for higher $s$:
even though the terms $s< n-1$ on the left hand side could be reduced lowering $n$ via \eqref{eq.lown},
the two terms $I^-_{s,1,s}$ and $I^-_{s+1,1,s}$ remain entangled.

Sum rules $\sum_{s=1}^n \alpha_{n,s}I^-_{n,j,s}$, $j>1$,  are derived from \eqref{eq.sumalpha}
by iterated differentiation of both sides w.r.t.\ $b$. The first of these steps is e.g.
\begin{equation}
\sum_{s=1}^n \alpha_{n,s}\hat I^-_{n,2,s}
= -\frac{\pi(n-1)!}{4(b\mu)^3}
-\frac{(-)^n(n-1)}{4}
\frac{(b\mu)^{n-3}}{(2\pi)^{n-1}}
\psi^{(n-1)}(\frac{b\mu}{2\pi})
-\frac{(-)^n}{4}
\frac{(b\mu)^{n-2}}{(2\pi)^n}
\psi^{(n)}(\frac{b\mu}{2\pi})
.
\end{equation}

\section{Partial Integration}\label{sec.pf}
Partial integration of $I^-_{3,j,s}$ may factorize the kernel fpr integrating $x/(x^2+b^2)^j$ 
and for deriving $x^2/(e^{\mu x}-1)^s$:
\begin{multline}
\int_0^\infty \frac{x^3}{(x^2+b^2)^j(e^{\mu x}-1)^s}dx
=
\frac{1}{2(1-j)}\frac{1}{(x^2+b^2)^{j-1}}
\frac{x^2}{(e^{\mu x}-1)^s}\mid_{x=0}^\infty
\\
-
\int_0^\infty 
\frac{1}{2(1-j)}\frac{1}{(x^2+b^2)^{j-1}}
\left[
\frac{2x}{(e^{\mu x}-1)^s}
-\frac{s\mu x^2 e^{\mu x}}{(e^{\mu x}-1)^{s+1}}
\right]dx.
\end{multline}
To avoid singularity at the first term of the right hand side (r.h.s.), $s\le 2$. 
For s=1
\begin{multline}
I^-_{3,j,1} = 
-
\int_0^\infty 
\frac{1}{2(1-j)}\frac{1}{(x^2+b^2)^{j-1}}
\left[
-\frac{x(\mu x-2)}{e^{\mu x}-1}
-\mu \frac{x^2}{(e^{\mu x}-1)^2}
\right]dx
\\
=
\frac{1}{2(j-1)}
\left[
2I^-_{1,j-1,1}
-\mu I^-_{2,j-1,1}
-\mu I^-_{2,j-1,2}
\right]
.
\label{eq.3j1}
\end{multline}
This allows to compute $I^-_{2,j,2}$-values based on $s=1$-values covered by sections \ref{sec.binom} and \ref{sec.1j1};
it may be simpler to compute these by iterated differentiation of \eqref{eq.dmu2} w.r.t.\ $b$.
\begin{rem}
For $s=2$
\begin{multline}
I^-_{3,j,2}
=
-\frac{1}{2(1-j)}\frac{1}{\mu^2 (b^2)^{j-1}}
-
\int_0^\infty 
\frac{1}{2(1-j)}\frac{1}{(x^2+b^2)^{j-1}}
\left[
\frac{2x}{(e^{\mu x}-1)^2}
-\frac{2\mu x^2 e^{\mu x}}{(e^{\mu x}-1)^3}
\right]dx
\\
=
\frac{1}{2(j-1)}\frac{1}{\mu^2 (b^2)^{j-1}}
+\frac{1}{j-1}
\int_0^\infty 
\frac{1}{(x^2+b^2)^{j-1}}
\left[
\frac{x}{(e^{\mu x}-1)^2}
-\frac{\mu x^2 e^{\mu x}}{(e^{\mu x}-1)^3}
\right]dx
.
\end{multline}
A new aspect arises as far as the two terms in the square brackets
do not converge individually, so a write-up as an infinite series is needed
for the expansion:
\begin{multline}
(j-1)I^-_{3,j,2}
=
\frac{1}{2}\frac{1}{\mu^2 (b^2)^{j-1}}
+
\int_0^\infty 
\frac{1}{(x^2+b^2)^{j-1}}
\frac{x}{e^{\mu x}-1}
[
\frac{1}{e^{\mu x}-1}
-\frac{\mu x }{e^{\mu x}-1}
-\frac{\mu x }{(e^{\mu x}-1)^2}
]dx
\\
=
\frac{1}{2}\frac{1}{\mu^2 (b^2)^{j-1}}
-\mu
I^-_{2,j-1,2}
+
\int_0^\infty 
\frac{1}{(x^2+b^2)^{j-1}}
\frac{x}{e^{\mu x}-1}
[
\frac{1}{e^{\mu x}-1}
-\frac{\mu x }{(e^{\mu x}-1)^2}
]dx
\\
=
\frac{1}{2}\frac{1}{\mu^2 (b^2)^{j-1}}
-\mu
I^-_{2,j-1,2}
+
\sum_{i=2}^\infty \frac{1}{i!}
\int_0^\infty 
\frac{1}{(x^2+b^2)^{j-1}}
\frac{x}{e^{\mu x}-1}
\frac{(\mu x)^i}
{(e^{\mu x}-1)^2}
dx
\\
=
\frac{1}{2}\frac{1}{\mu^2 b^{2j-2}}
-\mu
I^-_{2,j-1,2}
+
\sum_{i\ge 2} \frac{\mu^i}{i!}
I^-_{i+1,j-1,3}
.
\end{multline}
The same rephrased with the normalized integrals \eqref{def.norm}
\begin{equation}
(j-1)\hat I^-_{3,j,2}
=
\frac{1}{2}\frac{1}{(b\mu)^{2j-2}}
-
\hat I^-_{2,j-1,2}
+
\sum_{i=2}^\infty \frac{1}{i!}
\hat I^-_{i+1,j-1,3}
.
\end{equation}
\end{rem}

The partial fraction decomposition of $I^-_{4,j,s}$ may factorize the kernel as $x/(x^2+b^2)^j$ and $x^3/(e^{\mu x}-1)^s$:
\begin{multline}
I^-_{4,j,s} = -\frac{1}{2(j-1)}\frac{1}{(x^2+b^2)^{j-1}}\frac{x^3}{(e^{\mu x}-1)^s}\mid_{x=0}^\infty
\\
+
\int_0^\infty
\frac{1}{2(j-1)}\frac{1}{(x^2+b^2)^{j-1}}
[
\frac{3x^2}{(e^{\mu x}-1)^s}
-s\mu \frac{x^3 e^{\mu x}}{(e^{\mu x}-1)^{s+1}}
]
dx.
\end{multline}
In particular for $s=1$
\begin{multline}
I^-_{4,j,1} 
=
\int_0^\infty
\frac{1}{2(j-1)}\frac{1}{(x^2+b^2)^{j-1}}
[
-\frac{x^2(\mu x-3)}{e^{\mu x}-1}
-\mu \frac{x^3}{(e^{\mu x}-1)^2}
]
dx
\\
=
\frac{1}{2(j-1)}
[
-\mu I^-_{3,j-1,2}
+3 I^-_{2,j-1,1}
-\mu I^-_{3,j-1,1}
]
dx,
\end{multline}
basically equivalent to \eqref{eq.3j1}.
For $s=2$
\begin{multline}
I^-_{4,j,2} 
=
\int_0^\infty
\frac{1}{2(j-1)}\frac{1}{(x^2+b^2)^{j-1}}
[
-\frac{x^2(2\mu x-3)}{(e^{\mu x}-1)^2}
-2\mu \frac{x^3}{(e^{\mu x}-1)^3}
]
dx
\\
=
\frac{1}{2(j-1)}
[
-2\mu I^-_{3,j-1,2}
+3 I^-_{2,j-1,2}
-2\mu I^-_{3,j-1,3}
]
dx.
\end{multline}
This method of reducing terms with $s=3$ to terms $s\le 2$ 
compared 
to \eqref{eq.dmu3} 
has the disadvantage 
of needing more terms, but the advantage of covering general exponents $j$.

The rule following from these partial integrations that reduces $s=j=n-1$ to 3 terms with smaller $s$
turns out to be in general
\begin{equation}
I^-_{n,j,n-2}
=
\frac{1}{2(j-1)}
\left[
-(n-2)\mu I^-_{n-1,j-1,n-2}
+(n-1) I^-_{n-2,j-1,n-2}
-(n-2)\mu I^-_{n-1,j-1,n-1}
\right]
.
\end{equation}
The rule that reduces $s=j-1=n-2$ to 3 terms with smaller $s$
turns out to be in general
\begin{equation}
I^-_{n,j,n-3}
=
\frac{1}{2(j-1)}
\left[
-(n-3)\mu I^-_{n-1,j-1,n-3}
+(n-1) I^-_{n-2,j-1,n-3}
-(n-3)\mu I^-_{n-1,j-1,n-2}
.
\right]
.
\end{equation}

\section{Summary}
The integrals $I^\pm _{1,j,1}(b,\mu)$ are reduced to polygamma functions by partial differentiation
with respect to $b$ as outlined in Section \ref{sec.diffb} and
\ref{sec.1j1}.

The integrals $I^-_{n,j,s}(b,\mu)$ are reduced to either $I^-_{s,j,s}(b,\mu)$ or $I^-_{s+1,j,s}(b,\mu)$
with the binomial reduction of Section \ref{sec.binom}.

Parameters $s>1$ in $I^-_{s+1,j,s}(b,\mu)$ or $I^-_{s,j,s}(b,\mu)$
can be lowered 
by placing them into the last terms of right hand sides of
the last two equations
in Section \ref{sec.pf}.

\appendix
\section{Partial Integration: Canceling Divergences}
Partial integration is a generic strategy to generate recurrences of integrals.
The basic attempt is to factorize the kernel of $I^-_{1,j,s}$ by integrating the 
rational polynomial and by differentiating $1/(e^{\mu x}-1)^s$:
\begin{multline}
\int \frac{x}{(x^2+b^2)^j(e^{\mu x}-1)^s}dx
=
\frac{1}{2(1-j)}\frac{1}{(x^2+b^2)^{j-1}} \frac{1}{(e^{\mu x}-1)^s}
\\
+
\int \frac{1}{2(1-j)}\frac{1}{(x^2+b^2)^{j-1}} s\mu\frac{e^{\mu x}}{(e^{\mu x}-1)^{1+s}}dx,
\label{eq.div1}
\end{multline}
($s\le 1$ to converge on the l.h.s.).
Likewise for $I^-_{3,js}$ and $s\le 3$
\begin{multline}
\int \frac{x^3}{(x^2+b^2)^j}\frac{dx}{(e^{\mu x}-1)^s}
=
\frac{1}{2(1-j)(j-2)}\frac{b^2+(j-1)x^2}{(x^2+b^2)^{j-1}(e^{\mu x}-1)^s}
\\
+
\int \frac{1}{2(1-j)(j-2)}\frac{b^2+(j-1)x^2}{(x^2+b^2)^{j-1}} s\mu\frac{e^{\mu x}}{(e^{\mu x}-1)^{1+s}}dx.
\label{eq.div3}
\end{multline}
The prototypical problem for the $I^\pm$ integral family is that the first term on the r.h.s.\
becomes infinite at the lower limit $x\to 0^+$ and the new integral on the r.h.s.\ 
diverges reciprocally.
By plugging a dexterous linear combination of two different powers of $x$ in the
numerator these divergences cancel:
\begin{multline}
\leadsto
\int_0^\infty \frac{b^2x-(j-2)x^3}{(x^2+b^2)^j(e^{\mu x}-1)^s}dx
=
\frac{x^2/2}{(x^2+b^2)^{j-1}(e^{\mu x}-1)^s}\mid_0^\infty
+
\mu s
\int_0^\infty \frac{x^2/2 e^{\mu x}}{(x^2+b^2)^{j-1}(e^{\mu x}-1)^{s+1}}
\\
=
\frac{\mu s}{2}
\int_0^\infty \frac{x^2 e^{\mu x}}{(x^2+b^2)^{j-1}(e^{\mu x}-1)^{s+1}} dx,
\end{multline}
$s\le 1$ on the r.h.s.\ to converge.
With $s=1$, with the aid of \eqref{eq.lown} on the l.h.s.\ and deployment equivalent
to \eqref{eq.powexpmu} on the r.h.s.\ this means
\begin{equation}
(2-j)I^-_{1,j-1,1} +(j-1)b^2 I^-_{1,j,1}
=
\frac{\mu}{2}
(I^-_{2,j-1,1}+I^-_{2,j-1,2})
.
\label{eq.mix13}
\end{equation}
This is no new information because the same is obtained by
repeated differentiation of \eqref{eq.dmu2} w.r.t.\ $b$; note that $\alpha_{1,2}=\alpha_{2,2}$.

For $I^-_{5,j,s}$ the diverging variant is
\begin{multline}
\int \frac{x^5}{(x^2+b^2)^j}\frac{dx}{(e^{\mu x}-1)^s}
=
\frac{1}{2(1-j)(j-2)(j-3)}\frac{2b^4+\cdots x^2+\cdots x^4}{(x^2+b^2)^{j-1}(e^{\mu x}-1)^s}
\\
+
\int \frac{1}{2(1-j)(j-2)(j-3)}\frac{2b^4+\cdots x^2+\cdots x^4}{(x^2+b^2)^{j-1}} s\mu\frac{e^{\mu x}}{(e^{\mu x}-1)^{1+s}}dx.
\label{eq.div5}
\end{multline}
Mixing this with \eqref{eq.div3} yields
\begin{equation}
\leadsto
\int_0^\infty \frac{2b^2x^3-(j-3)x^5}{(x^2+b^2)^j(e^{\mu x}-1)^s}dx
=
\frac{\mu s}{2}
\int_0^\infty \frac{x^4 e^{\mu x}}{(x^2+b^2)^{j-1}(e^{\mu x}-1)^{s+1}} dx,
\end{equation}
which requires $s\le 3$ to converge. 
\begin{equation}
(3-j)I^-_{3,j-1,s}
+(j-1)b^2 I^-_{3,j,s}
=\frac{\mu s}{2}(I^-_{4,j-1,s}+I^-_{4,j-1,s+1}),\quad s\le 3
\label{eq.mix35}
\end{equation}
\begin{rem}
For $j=1$, $s=1$, the l.h.s.\ is easily evaluated with \eqref{eq.noden}
and both sides reduce to a well-known integral \cite[3.423.2]{GR}.
\end{rem}

\begin{rem}
Mixing \eqref{eq.div5} and \eqref{eq.div1} (for $s\le 1$ to converge)
\begin{multline}
\int_0^\infty \frac{2b^4x-(j-2)(j-3)x^5}{(x^2+b^2)^j(e^{\mu x}-1)^s}dx
=
\frac{\mu s}{2}\int_0^\infty \frac{x^2[2b^2+(j-2)x^2]}{(x^2+b^2)^{j-1}} \frac{e^{\mu x}}{(e^{\mu x}-1)^{s+1}}dx
\end{multline}
yields more convoluted correlations:
\begin{multline}
(-j^2+4j-3)I^-_{3,j-1,s}
+(j^2b^2-4jb^2+3b^2)I^-_{1,j-1,s}
-b^4(1+j^2-4j)I^-_{1,j,s}
\\
=
\frac{\mu s}{2}\left[
2b^2(I^-_{2,j-1,s}+I^-_{2,j-1,s+1})
+(j-2)(I^-_{4,j-1,s}+I^-_{4,j-1,s+1})
\right]
.
\label{eq.mix15}
\end{multline}
This is no new information but a mere linear combination of \eqref{eq.mix13} and \eqref{eq.mix35}.
\end{rem}
\section{Numerical Integration}

There are parameter sets where the previous chapters 
cannot reduce the integrals to polygamma-functions. In these
cases a standard numerical approach
substitutes 
$z=x/(x+c)$, 
$x=cz/(1-z)$
which maps the infinite $x$-interval to the unit $z$-interval:
\begin{equation}
I^\pm _{n,j,s}(b,\mu) 
=
\frac{c^{n+1}}{(c^2+b^2)^j} \int_0^1 \frac{z^n}{(1-z)^{n+2-2j}(z^2+\frac{b^2}{c^2+b^2}(1-2z))^j(e^{\mu cz/(1-z)}\pm 1)^s}dz
.
\label{eq.zunit}
\end{equation}
$c$ is a free tuning parameter that according to some homespun experimentation 
ought to be roughly $\sqrt{b\mu}$ to
shift the maximum of the integral kernel into the mid-region $z\approx 0.5$.

\begin{table}
\begin{tabular}{lllrrr}
$n$ & $j$ & $s$ & $b\mu$ & $I^-_{n,j,s}(b\mu,1)$ & $I^+_{n,j,s}(b\mu,1)$\\
\hline
\input{takec1.out}
\end{tabular}
\caption{Reference values of $I^\pm_{n,j,s}(1.,1.)$}
\end{table}

\begin{table}
\begin{tabular}{lllrrr}
$n$ & $j$ & $s$ & $b\mu$ & $I^-_{n,j,s}(b\mu,1)$ & $I^+_{n,j,s}(b\mu,1)$\\
\hline
\input{takecPi.out}
\end{tabular}
\caption{Reference values of $I^\pm_{n,j,s}(\pi,1.)$}
\end{table}

\begin{table}
\begin{tabular}{lllrrr}
$n$ & $j$ & $s$ & $b\mu$ & $I^-_{n,j,s}(b\mu,1)$ & $I^+_{n,j,s}(b\mu,1)$\\
\hline
\input{takec2Pi.out}
\end{tabular}
\caption{Reference values of $I^\pm_{n,j,s}(2\pi,1.)$}
\end{table}

The two files of the C-source code in the ancillary directory
implement \eqref{eq.zunit} 
with a Boole-rule; the are compiled with
\begin{verbatim}
gcc -O2 -o digamma_assoc -D TEST [-D DEBUG] digamma_assoc.c
\end{verbatim}
to create the executable \verb+digamma_assoc+ . A \verb+main+
program is included if the \verb+TEST+ C-preprocessor (CPP) option is used,
and traces of the number of sampling points in $z$ and $I^-$
are printed if the \verb+DEBUG+ CPP option is used.

The program is called with
\begin{verbatim}
digamma_assoc [-p] [-n n] [-j j] [-s s] [-m mu] [-b b] [-e epsilon]
\end{verbatim}
where the optional command line arguments are $n$, $j$, $s$, $\mu$
and $b$ (defaults 1, 1, 1, 1 and $2\pi$) and \texttt{epsilon} is a requested absolute error (default $10^{-10}$) of the $I^\pm$
value.
If the \texttt{-p} option is used, $I^+_{n,j,s}(b,\mu)$ is calculated, otherwise
$I^-_{n,j,s}(b,\mu)$.

\bibliographystyle{amsplain}
\bibliography{all}

\end{document}